\documentclass[10pt]{article}

\usepackage{pifont}
\usepackage{amsmath}
\usepackage{cases}
\usepackage{amsfonts}
\usepackage{latexsym}
\usepackage{amssymb}
\usepackage{stmaryrd}
\usepackage{overpic}
\usepackage{color}
\usepackage{abstract}
\usepackage{indentfirst}
\usepackage{booktabs}
\usepackage{url}
\usepackage{verbatim}
\usepackage{diagbox}
\usepackage{multirow} 
\usepackage{array} 
\usepackage{graphicx}  
\usepackage{float}  
\usepackage{subfigure}

\newcommand{\btau}{\boldsymbol{\tau}}

\newcommand{\bvarepsilon}{\boldsymbol{\varepsilon}}
\newcommand{\bsigma}{\boldsymbol{\sigma}}

\newcommand{\bzero}{\boldsymbol{0}}
\newcommand{\bv}{\boldsymbol{v}}
\newcommand{\bx}{\boldsymbol{x}}
\newcommand{\fb}{\boldsymbol{f}}
\newcommand{\bu}{\boldsymbol{u}}
\newcommand{\bH}{\boldsymbol{H}}
\newcommand{\bQ}{\boldsymbol{Q}}

\newcommand{\bg}{\boldsymbol{g}}
\newcommand{\bn}{\boldsymbol{n}}

\newcommand{\ds}{\mathrm{d}s}

\DeclareFontEncoding{FMS}{}{}
\DeclareFontSubstitution{FMS}{futm}{m}{n}
\DeclareFontEncoding{FMX}{}{}
\DeclareFontSubstitution{FMX}{futm}{m}{n}
\DeclareSymbolFont{fouriersymbols}{FMS}{futm}{m}{n}
\DeclareSymbolFont{fourierlargesymbols}{FMX}{futm}{m}{n}
\DeclareMathDelimiter{\VERT}{\mathord}{fouriersymbols}{152}{fourierlargesymbols}{147}

\newtheorem{theorem}{Theorem}[section]

\newtheorem{example}[theorem]{Example}

\newtheorem{remark}[theorem]{Remark}

\begin{document}
	
	\begin{center}
		\large\bf Randomized Neural Networks with Petrov-Galerkin Methods\\ for Solving Linear Elasticity Problems
	\end{center}
	
	\begin{center}
		{\large\sc Yong Shang}\footnote{School of Mathematics and Statistics, Xi'an Jiaotong University, Xi'an,
			Shaanxi 710049, P.R. China. E-mail: {\tt fsy2503.xjtu@xjtu.edu.cn}},\quad
		{\large\sc Fei Wang}\footnote{School of Mathematics and Statistics, Xi'an Jiaotong University,
			Xi'an, Shaanxi 710049, China. The work of this author was partially supported by
			the National Natural Science Foundation of China (Grant No.\ 12171383). Email: {\tt feiwang.xjtu@xjtu.edu.cn}}
	\end{center}

	\begin{quote} 
		\noindent{}{\bf Abstract}: We develop the Randomized Neural Networks with Petrov-Galerkin Methods (RNN-PG methods) to solve linear elasticity problems. RNN-PG methods use Petrov-Galerkin variational framework, where the solution is approximated by randomized neural networks and the test functions are piecewise polynomials. Unlike conventional neural networks, the parameters of the hidden layers of the randomized neural networks are fixed randomly, while the parameters of the output layer are determined by the least square method, which can effectively approximate the solution. We also develop mixed RNN-PG methods for linear elasticity problems, which ensure the symmetry of the stress tensor and avoid locking effects. We compare RNN-PG methods with the finite element method, the mixed discontinuous Galerkin method, and the physics-informed neural network on several examples, and the numerical results demonstrate that RNN-PG methods achieve higher accuracy and efficiency.
		

		{\bf Keywords:} Petrov-Galerkin formulation, randomized neural networks, linear elasticity problem, least-square method.
	\end{quote}
	
	\section{Introduction}
	
	Linear elasticity problems deal with the displacement, strain, and stress fields in a solid under external forces and boundary conditions. These problems are fundamental for engineering because they help design and analyze structures and materials that can withstand loads without excessive deformation or fracture. Computing linear elasticity problems accurately is crucial for many applications. However, the low-order finite element method (FEM) may encounter a ``locking” issue for nearly incompressible materials (\cite{Babuska1992locking}), which prevents the method from capturing the correct deformation modes of the structure. Several “locking-free” approaches have been developed, such as nonconforming methods (\cite{Kouhia1995nonconforming}), discontinuous Galerkin (DG) method (\cite{Hansbo2002dg}), mixed DG method (\cite{Wang2020mixdg}) and others. These methods tend to increase the number of degrees of freedom and become computationally challenging in higher dimensions.


	Artificial neural network methods have a strong approximation ability, which offers a new approach to solving partial differential equations (PDEs). In particular, they show remarkable advantages for solving high-dimensional problems. The application of neural networks to solve differential equations was pioneered by \cite{Lagaris1998ann}. In recent years, deep neural networks-based methods for solving PDEs have been developed rapidly, such as the Deep Ritz Method (\cite{E2018drm}), the Deep Galerkin Method (\cite{sirignano2018dgm}), and the Physics-Informed Neural Networks (PINNs, \cite{Raissi2019pinn}). These NN-based methods train neural networks by minimizing a loss function that incorporates the underlying PDEs and initial/boundary conditions. PINNs have been applied to various problems in fluid mechanics (\cite{Sun2020pcdl,Jin2021NSFnets,Mao2021NSFnets}), solid mechanics (\cite{Haghighat2020ex2solid,Goswami2021phasefield,Rao2021elastodynamics}) and other fields (\cite{Kadeethum2020diffusivity,Cai2021heattransfer}). Some research works on PINNs for linear elasticity problems are \cite{Guo2022ex2elasticity,Rezaei2022mixpinn,Zhang2022internal,Roy2023padl,Samaniego2023energy,Vahab2023biharmonic}.

	%

	These DNN-based methods have achieved significant advancements in solving high-dimensional and irregular domain problems. However, the training process can be computationally expensive, and it is challenging to reach the optimal state of the problem with the existing training algorithms. To overcome this challenge, randomized neural networks-based methods have been proposed for solving PDEs, such as Local Extreme Learning Machines and Domain Decomposition Method (\cite{dong2021elm}), Randomized Neural Networks with Petrov-Galerkin Methods (RNN-PG methods, \cite{Shang2022dpgm}), Local Randomized Neural Networks with Discontinuous Galerkin Methods (\cite{Sun2022lrrdg}), and Random Feature Method (\cite{chen2022randomfm}). In these methods, the unknown variables are approximated by randomized neural networks, and parameters are solved by least-square computation rather than training algorithms, so they can approximate the solutions more accurately with less computational cost.

	In this paper, we develop Randomized Neural Networks with Petrov-Galerkin Methods to solve linear elasticity problems. The paper is organized as follows. Section \ref{sec:RNN-PG} introduces the RNN-PG method for solving a linear elasticity problem. In order to present a scheme for preserving the symmetry of the stress tensor, we introduce Mixed RNN-PG methods based on different mixed formulations in Section 3. Numerical examples are shown in Section 4 with a comparison of FEM, Mixed DG method, and PINNs. The last section concludes the paper with some discussions.

	

	\section{RNN-PG method for a linear elasticity problem}
	\label{sec:RNN-PG}
	
	In this section, we introduce Randomized Neural Networks with Petrov-Galerkin (RNN-PG) Methods for solving linear elasticity problems.
	
	Let $\Omega$ be a bounded domain in $\mathbb{R}^d$ ($d=2,3$) with boundary $\Gamma = \Gamma_D \cup \Gamma_N $ and ${\Gamma}_D \cap {\Gamma}_N = \emptyset$. Consider a linear elasticity problem of the form:
	\begin{numcases}{}
		\mathcal{A}\bsigma - \bvarepsilon(\bu)=0\;\;\;\;\qquad\;\;\;\;\;\; \text{in}\; \Omega\; , \label{prb1}\\
		- {\rm div}\, \bsigma = \fb  \;\;\;  \qquad \qquad\quad\;\,\text{in} \; \Omega\;  ,\label{prb2}\\
		\bu= {\bg_D}\qquad\;\;\;\,\; \qquad \qquad\;\,\;\text{on}\; \Gamma_D \; ,\label{bc1}\\
		\bsigma \bn = {\bg_N} \qquad \qquad \qquad \,\;\;\;\;\text{on}\; \Gamma_N, \label{bc2}\;
	\end{numcases}
	where displacemnt $\bu:\Omega \rightarrow \mathbb{R}^d$ and the stree tensor $\bsigma : \Omega \rightarrow \mathbb{S}$ with $\mathbb{S}$ denoting the space of real symmetric matrices of order $d \times d$.	Here, $\bvarepsilon(\bu) =(\nabla \bu + (\nabla \bu)^t)/2$ is the strain tensor, $\lambda$ and $\mu$ are Lam\'e coefficients,  and compliance tensor $	\mathcal{A}:\mathbb{S}\rightarrow \mathbb{S}$ is given by 
	\begin{equation}
		\mathcal{A}\bsigma = \frac{1}{2\mu}\left(\bsigma - \frac{\lambda}{2\mu + d\lambda}\text{tr}(\bsigma)I_d\right).
	\end{equation}
	
	The weak formulation of the linear elasticity problem \eqref{prb1}--\eqref{bc2} is: Find $\bu\in H^1_{D,g_D}(\Omega;\mathbb{R}^d) =  \{u\in H^1(\Omega;\mathbb{R}^d); \;u|_{\Gamma_D} = g_D\}$ such that
	\begin{align}
		a(\bu,\bv) &= l(\bv) \qquad \;\forall \bv  \in H^1_{D,0}(\Omega;\mathbb{R}^d) = \{v\in H^1(\Omega;\mathbb{R}^d);\; v|_{\Gamma_D} = 0\},\label{weak}
	\end{align}
	where
	\begin{align*}
		a(\bu,\bv)&=\int_\Omega \big(2\mu  \bvarepsilon(\bu)  :  \bvarepsilon(\bv) + \lambda \nabla \cdot \bu \,\nabla \cdot \bv \big)	\, {\rm d}\bx,\\
		l(\bv) &=\int_\Omega \fb\cdot \bv\,{\rm d}\bx\,+ \int_{\Gamma_N} \bg_N\cdot \bv\,\ds.
	\end{align*}

	Instead of using piecewise polynomial functions to approximate the solution of the aforementioned problem, we employ neural networks to approximate the displacement vector $\bu$. To solve the problem \eqref{weak}, we utilize the RNN-PG method and aim to find a neural network $\bu_\rho \in U_\rho $ such that
	\begin{align}
		a(\bu_\rho,\bv) &= l(\bv) \qquad \;\;\;\forall \bv  \in V_h \subset H^1_{D,0}(\Omega;\mathbb{R}^d), \label{prb3} \\
		\bu_{\rho}(\bx_k) &=  \bg_D(\bx_k)\;\;\,\,\;\;{\text{for some points}}\ \bx_k\in  \Gamma_D,\; k=1,2,\cdots,N_b. \label{bc3}\ 
	\end{align}
	Here, equation \eqref{bc3} is used to enforce the Dirichlet boundary condition \eqref{bc1} by selecting random samples $\{\bx_k\}_{k=1}^{N_b}$ according to the uniform distribution $\mathcal{U}(\Gamma_D)$. In addition, $V_h$ can be any finite-dimensional function space that effectively approximates $H^1_{D,0}(\Omega;\mathbb{R}^d)$, while $U_\rho $ consists of randomized neural networks $\mathbf{\Phi} \in H^1(\Omega;\mathbb{R}^d)$ with fixed parameters for input layer and hidden layers, and adjustable parameters for the output layer. The concept of randomized neural networks was initially introduced in \cite{Pao1994RNN4}, and extreme learning machine (ELM) is one example (\cite{Huang2006ELMtheorandapp}). In comparison to feedforward neural networks, ELM preserves generalization capability when appropriate activation functions are employed and suitable initialization methods for fixed parameters are used (\cite{Liu2014ELMfeasible}).

	\begin{figure}[!htbp] 		
		\centering
		\includegraphics[scale=0.5]{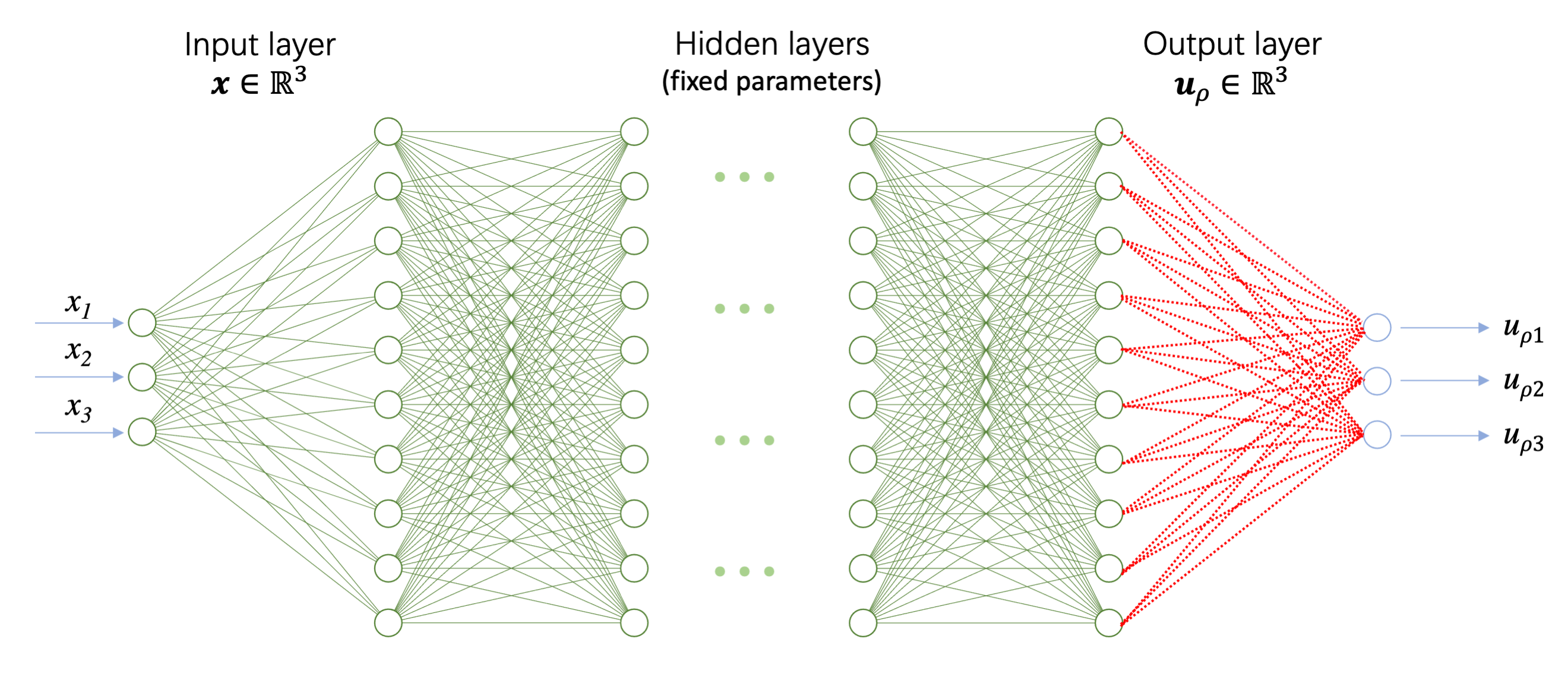}
		\caption{Network structure of $\mathbf{\Phi}:\mathbb{R}^{3}\rightarrow \mathbb{R}^3$, the solid green line represents the parameters of the neural network that are randomly initialized and fixed thereafter, while the dotted red line refers to the adjustable parameters. }
		\label{figure1}
	\end{figure}
	
	Let $D$ be the depth and $\rho$ be the activation function. A fully connected feedforward neural network $\mathbf{\Phi}$ is defined by
	\begin{align*}
		&\mathbf{\Phi}_0(\bx) = \bx,\nonumber\\
		&\mathbf{\Phi}_{l}(\bx) = \rho(\mathbf{W}_l\mathbf{\Phi}_{l-1}+\mathbf{b}_l) \qquad \text{for}\ l=1,\cdots,D-1,\label{neural}\\
		&\mathbf{\Phi}:= \mathbf{\Phi}_D(\bx) =\mathbf{W}_D\mathbf{\Phi}_{D-1},\nonumber
	\end{align*}
	where $ \left\{ \mathbf{W}_l = \left(w_{ij}^{(l)} \right)\in \mathbb{R}^{n_l \times n_{l-1}}, \mathbf{b}_l = \left(b_i^{(l)}\right)\in \mathbb{R}^{n_l} \right\}$ are the weight parameters initialized with a uniform distribution $\mathcal{U}(-r, r)$ in the $l$-th layer, and $r\in \mathbb{R}$. 
	
	When $d=3$, let $\bu_\rho$ be the output of $\mathbf{\Phi}$ with 3 input neurons and 3 output neurons shown in Figure \ref{figure1}. Denote $\Phi_j^u :=\mathbf{\Phi}_{D-1}^j$, $u_{1\rho}^j :=w_{1j}^{D}$ , $u_{2\rho}^j :=w_{2j}^{D}$  and $u_{3\rho}^j :=w_{3j}^{D}$ for $j = 1,\cdots,n_{D-1}$. Then 
	\begin{equation*}
		\bu_\rho(\bx) = \mathbf{W}_D\mathbf{\Phi}_{D-1} =
		\left(\sum\limits_{j=1}^{ n_{{{D}}-1}}u_{1\rho}^{j} \Phi_{j}^{u}(\bx), \sum\limits_{j=1}^{ n_{{{D}}-1}}u_{2\rho}^{j} \Phi_{j}^{u}(\bx),\sum\limits_{j=1}^{ n_{{{D}}-1}}u_{3\rho}^{j} \Phi_{j}^{u}(\bx)\right)^t := (u_{\rho1},u_{\rho2},u_{\rho3})^t,
	\end{equation*}
	where the variables $u_{1\rho}^j$, $u_{2\rho}^j$, and $u_{3\rho}^j$ are to be determined for $j = 1, \ldots, n_{D-1}$.
	
	Test functions can be chosen from $\{\bv_i\}_{i=1}^{N_h}  \in V_h \subset H^1_{D,0}(\Omega;\mathbb{R}^d)$. 
	Therefore, problem \eqref{prb3}-\eqref{bc3} becomes: Find $u_{1\rho}^j,u_{2\rho}^j$ and $u_{3\rho}^j$, $j = 1,\cdots,n_{D-1}$, such that
	\begin{align}
		a\left(\left(\sum\limits_{j=1}^{ n_{{{D}}-1}}u_{1\rho}^{j} \Phi_{j}^{u}(\bx),	\sum\limits_{j=1}^{ n_{{{D}}-1}}u_{2\rho}^{j} \Phi_{j}^{u}(\bx), \sum\limits_{j=1}^{ n_{{{D}}-1}}u_{3\rho}^{j} \Phi_{j}^{u}(\bx)\right)^t,\bv_i\right) &= l(\bv_i) \quad\;\;\, \text{for} \  i = 1,\cdots,N_h,  \label{linear1}\\
		\left(\sum\limits_{j=1}^{ n_{{{D}}-1}}u_{1\rho}^{j} \Phi_{j}^{u}(\bx_k),	\sum\limits_{j=1}^{ n_{{{D}}-1}}u_{2\rho}^{j} \Phi_{j}^{u}(\bx_k), \sum\limits_{j=1}^{ n_{{{D}}-1}}u_{3\rho}^{j} \Phi_{j}^{u}(\bx_k)\right)^t &= \bg_D(\bx_k)  \;\; \text{for} \  k = 1,\cdots,N_b. \label{linear2}\end{align}
	
	Finally, we can obtain the solution $\bu_\rho$ by solving a least-squares problem with the corresponding linear system \eqref{linear1}-\eqref{linear2}.

	\section{Mixed RNN-PG methods for the linear elasticity problem}
	
	Elasticity stress measures how much a material can deform under external forces without breaking or losing its original shape. This is an important concept for engineers who design structures such as buildings and aircraft. Mixed methods can approximate both displacement and stress tensor at the same time. However, the stress tensor has to be symmetric according to the conservation of angular momentum principle (\cite{Arnold2008stress}). Designing stable mixed finite element methods with symmetric stress is a challenging task for linear elasticity problems, and many research works have investigated this topic, see  \cite{Arnold2007mixfem,Arnold2008fem,Cockburn2010element,Hu2015higher,Hu2015symmix,Hu2016lower,Wu2017ipmix,Wang2020mixdg} and the references therein.
	
	In this section, we develop the Mixed RNN-PG (M-RNN-PG) methods, which use separate neural networks to approximate displacement $\bu$ and stress tensor $\bsigma$, respectively, and maintain the symmetric property automatically.

	
	For example, when $d=2$, the unkown varaibles $\bu$ and $\bsigma$ are approxiamted by $\bu_\rho$ and $\bsigma_\rho$  as following.
	Denote $\bu  = \begin{pmatrix}  u_1\\u_2 \end{pmatrix}$, $\bsigma = 	\begin{pmatrix}   \sigma_{11}& \sigma_{12}\\\sigma_{21} &\sigma_{22} \end{pmatrix}$, 
	then  $u_1$, $u_2$, $\sigma_{11}$, $\sigma_{12}$, $\sigma_{21}$, $\sigma_{22}$  were appximated by $ \boldsymbol{\Phi}^u_1 $, $ \boldsymbol{\Phi}^u_2$, $\boldsymbol{\Phi}^{\sigma}_1$, $\boldsymbol{\Phi}^{\sigma}_2$, $\boldsymbol{\Phi}^{\sigma}_2$, $\boldsymbol{\Phi}^{\sigma}_3$, separately. By using the same approximation of $\sigma_{12}$ and $\sigma_{21}$ by $\boldsymbol{\Phi}^{\sigma}_2$, the symmetry property of numerical solution $\bsigma_\rho$ is preserved automatically. As a similar setting in Section \ref{sec:RNN-PG}, the randomized neural networks are defined by
	\begin{equation*}
		\bu_\rho =\mathbf{W}_{\tilde{D}}\boldsymbol{\Phi}^u_{{\tilde{D}}-1} = 
		\left[\begin{array}{c}\sum\limits_{j=1}^{ n_{{\tilde{D}}-1}}u_{1\rho}^{j} \Phi_{j}^{u}(\bx) \\
			\sum\limits_{j=1}^{ n_{{\tilde{D}}-1}}u_{2\rho}^{j} \Phi_{j}^{u}(\bx)\end{array}\right] 
		:= \left[\begin{array}{c}  \boldsymbol{\Phi}^u_1 \\ \boldsymbol{\Phi}^u_2\end{array}\right],
	\end{equation*}
	
	\begin{equation*}
		\bsigma_\rho =\mathbf{W}_{\hat{D}}\boldsymbol{\Phi}^{\sigma}_{{\hat{D}}-1} = 
		\left[\begin{array}{c}\sum\limits_{j=1}^{ n_{{\tilde{D}}-1}}{\sigma}_{1\rho}^{j} \Phi_{j}^{{\sigma}}(\bx) \\
			\sum\limits_{j=1}^{ n_{{\hat{D}}-1}}{\sigma}_{2\rho}^{j} \Phi_{j}^{{\sigma}}(\bx) \\ \sum\limits_{j=1}^{ n_{{\hat{D}}-1}}{\sigma}_{3\rho}^{j} \Phi_{j}^{{\sigma}}(\bx) \end{array}\right]
		:= \left[\begin{array}{c}  \boldsymbol{\Phi}^{\sigma}_1 \\  \boldsymbol{\Phi}^{\sigma}_2  \\  \boldsymbol{\Phi}^{\sigma}_3 \end{array}\right] .
	\end{equation*}
	
	To obtain different mixed formulations, one can multiply test functions on both sides of the equations \eqref{prb1}--\eqref{prb2}, followed by integration by parts and the boundary conditions \eqref{bc1}--\eqref{bc2}, then four weak formulations are given as follows.
	
	{\bf Mixed Formulation 1}:  
	Find $(\bsigma,\bu) \in L^2(\Omega;\mathbb{S}) \times H^1_{D,g_D}(\Omega;\mathbb{R}^d)$ such that 
	\begin{align*}
		\int_{\Omega}\mathcal{A}\bsigma:\btau \ {\rm d}\bx - \int_{\Omega}\bvarepsilon(\bu):\btau\ {\rm d}\bx &= 0 \qquad\qquad\qquad\qquad\quad\qquad \;\;\quad\forall \btau\in L^2(\Omega;\mathbb{S}), \\
		\int_{\Omega} \bsigma:\bvarepsilon(\bv) \ {\rm d}\bx  &= \int_{\Omega} \fb \cdot\bv \ {\rm d}\bx +\int_{\Gamma_N}\bg_N\cdot \bv \ \ds \quad \forall \bv \in H^1_{D,0}(\Omega;\mathbb{R}^d).
	\end{align*}
	
	{\bf Mixed Formulation 2}:  
	Find $(\bsigma,\bu)  \in \bH_{N,g_N}(\mathrm{div},\Omega;\mathbb{S})  \times L^2(\Omega;\mathbb{R}^d)$ such that 
	\begin{align*}
		\int_{\Omega}\mathcal{A}\bsigma:\btau \ {\rm d}\bx + \int_{\Omega}\bu \cdot {\rm div} \btau \ {\rm d}\bx&= \int_{\Gamma_D}\bg_D\cdot \btau \bn\ \ds \quad \;\;\,\forall \btau \in \bH_{N,0}(\mathrm{div},\Omega;\mathbb{S}), \\
		-\int_{\Omega}  {\rm div} \bsigma \cdot v\ {\rm d}\bx &= \int_{\Omega} \fb \cdot \bv \ {\rm d}\bx \qquad\quad\;\; \forall \bv \in L^2(\Omega;\mathbb{R}^d).
	\end{align*}
	
	{\bf Mixed Formulation 3}:  Find $(\bsigma,\bu)  \in \bH_{N,g_N}(\mathrm{div},\Omega;\mathbb{S})  \times H^1_{D,g_D}(\Omega;\mathbb{R}^d)$ such that 
	\begin{align*}
		\int_{\Omega}\mathcal{A}\bsigma:\btau \ {\rm d}\bx - \int_{\Omega}\bvarepsilon(\bu):\btau\ {\rm d}\bx &= 0 \qquad\qquad\qquad\qquad \forall \btau\in L^2(\Omega;\mathbb{S}), \\
		-\int_{\Omega}  {\rm div} \bsigma \cdot v\ {\rm d}\bx &= \int_{\Omega} \fb \cdot \bv \ {\rm d}\bx \qquad\quad\;\; \,\forall \bv \in L^2(\Omega;\mathbb{R}^d).
	\end{align*}
	
	{\bf Mixed Formulation 4}: Find $(\bsigma,\bu)  \in L^2(\Omega;\mathbb{S}) \times L^2(\Omega;\mathbb{R}^d)$ such that
	\begin{align*}
		\int_{\Omega}\mathcal{A}\bsigma:\btau \ {\rm d}\bx + \int_{\Omega}\bu \cdot {\rm div} \btau \ {\rm d}\bx&= \int_{\Gamma_D}\bg_D\cdot \btau \bn\ \ds \qquad\qquad\qquad\;\, \forall \btau \in \bH_{N,0}(\mathrm{div},\Omega;\mathbb{S}), \\
		\int_{\Omega} \bsigma:\bvarepsilon(\bv) \ {\rm d}\bx  &= \int_{\Omega} \fb \cdot\bv \ {\rm d}\bx +\int_{\Gamma_N}\bg_N\cdot \bv \ \ds \quad \forall \bv \in H^1_{D,0}(\Omega;\mathbb{R}^d).
	\end{align*}
	Here,
	\begin{align*}
		H^1_{D,0}(\Omega;\mathbb{R}^d) &= \{\bv\in H^1(\Omega;\mathbb{R}^d);\; \bv|_{\Gamma_D} = 0\},\\ 
		H^1_{D,g_D} (\Omega;\mathbb{R}^d) &= \{\bv\in H^1(\Omega;\mathbb{R}^d);\; \bv|_{\Gamma_D} = g_D\},\\ 
		\bH_{N,0}(\mathrm{div},\Omega;\mathbb{S}) &= \{\btau\in \bH(\mathrm{div},\Omega;\mathbb{S});\; \langle\btau \bn ,\bv\rangle = 0\; \forall \bv \in H^1_{D,0}(\Omega;\mathbb{R}^d)\},\\
		\bH_{N,g_N}(\mathrm{div},\Omega;\mathbb{S})& = \{\btau\in \bH(\mathrm{div},\Omega;\mathbb{S});\; \langle\btau \bn ,\bv\rangle = \bg_N\; \forall \bv \in H^1_{D,0}(\Omega;\mathbb{R}^d)\},
	\end{align*}
	and $\langle\cdot ,\cdot\rangle$ denotes the duality between $H^{-1/2}(\Gamma;\mathbb{R}^d)$ and $H^{1/2}(\Gamma;\mathbb{R}^d)$.
	
	\begin{remark}
		The weak form of the Mixed Formulation 4 incorporates all the boundary conditions, so there is no need for additional treatment of the boundary. Moreover, the weak form eliminates the computation of derivatives, which may enhance the accuracy of the results.
	\end{remark}
	
	Let us take Mixed Formulation 4 as an example to illustrate how the M-RNN-PG method works. The equivalent form 
	reads: Find $(\bsigma,\bu)  \in L^2(\Omega;\mathbb{S}) \times L^2(\Omega;\mathbb{R}^d)$ such that
	\begin{align}
		\mathcal{L}\big((\bsigma,\bu);(\btau,\bv)\big) = \mathcal{F}(\btau,\bv)\qquad \forall (\btau,\bv)\in \bH_{N,0}(\mathrm{div},\Omega;\mathbb{S}) \times H^1_{D,0}(\Omega;\mathbb{R}^d).\label{prb4}
	\end{align}
	where 
	\begin{align*}
		\mathcal{L}\big((\bsigma,\bu);(\btau,\bv)\big) &= \int_{\Omega} \big(\mathcal{A}\bsigma:\btau+\bu \cdot {\rm div} \btau+ \bsigma:\bvarepsilon(\bv) \big) \ {\rm d}\bx, \\
		\mathcal{F}(\btau,v) &= \int_{\Omega} \fb \cdot\bv \ {\rm d}\bx +\int_{\Gamma_N}\bg_N\cdot \bv \ \ds + \int_{\Gamma_D}\bg_D\cdot \btau \bn\ \ds .
	\end{align*}
	
	Therefore, the M-RNN-PG method is: Find neural networks $\bu_\rho \in U_\rho$ and $\bsigma_\rho \in \bQ_\rho$ such that
	\begin{align}
		\mathcal{L}\big((\bsigma_\rho,\bu_\rho);(\btau_i,\bv_k)\big) = \mathcal{F}(\btau_i,\bv_k)\qquad \forall (\btau_i,\bv_k)\in \bQ_h\times V_h,\label{prb10}
	\end{align}
	where $\bQ_h$ and $V_h$ can be chosen as any proper finite-dimensional function spaces, $U_\rho$ and $\bQ_\rho$ consist of randomized neural networks $\bu_\rho \in L^2(\Omega;\mathbb{R}^d)$ and $\bsigma_\rho \in  L^2(\Omega;\mathbb{S})$, respectively.
	
	Thus, problem \eqref{prb10} becomes: Find $u_{1\rho}^j ,u_{2\rho}^j  $ with $j = 1,\cdots,n_{\tilde{D}-1}$, and $\sigma_{1\rho}^{j}, \sigma_{2\rho}^{j},\sigma_{3\rho}^{j}$ with $j = 1,\cdots,n_{{\hat{D}}-1}$ such that
	\begin{align}\label{prb11}
		\mathcal{L}\left(\left(
		\begin{pmatrix}   
			\sum\limits_{j=1}^{ n_{{\hat{D}}-1}}{\sigma}_{1\rho}^{j} \Phi_{j}^{{\sigma}}(\bx) & 
			\sum\limits_{j=1}^{ n_{{\hat{D}}-1}}{\sigma}_{2\rho}^{j} \Phi_{j}^{{\sigma}}(\bx) \\
			\sum\limits_{j=1}^{ n_{{\hat{D}}-1}}{\sigma}_{2\rho}^{j} \Phi_{j}^{{\sigma}}(\bx) & \sum\limits_{j=1}^{ n_{{\hat{D}}-1}}{\sigma}_{3\rho}^{j} \Phi_{j}^{{\sigma}}(\bx)
		\end{pmatrix},
		\begin{pmatrix} \sum\limits_{j=1}^{ n_{{\tilde{D}}-1}}u_{1\rho}^{j} \Phi_{j}^{u}(\bx) \\
			\sum\limits_{j=1}^{ n_{{\tilde{D}}-1}}u_{2\rho}^{j} \Phi_{j}^{u}(\bx)\end{pmatrix} 
		\right);(\btau_i,\bv_k)\right) = \mathcal{F}(\btau_i,\bv_k) \;\;\forall (\btau_i,\bv_k)\in \bQ_h\times V_h.
	\end{align}
	For test functions, denote $\btau_i = \begin{pmatrix}   \tau_{11}& \tau_{12}\\\tau_{12} &\tau_{22} \end{pmatrix}$ and $\bv_k = \begin{pmatrix}  v_1\\v_2 \end{pmatrix}$, then we can take $(\btau_i,\bv_k) \in \bQ_h\times V_h$ in the forms of 
	$$ 
	\left( \begin{pmatrix}  \tau_{11}& 0\\0&0 \end{pmatrix}, \begin{pmatrix}  v_1\\0 \end{pmatrix}\right), \quad \left( \begin{pmatrix}   \tau_{11}& 0\\0&0 \end{pmatrix}, \begin{pmatrix}  0\\v_2 \end{pmatrix}\right), \quad \left( \begin{pmatrix}   0& \tau_{12}\\\tau_{12}&0 \end{pmatrix}, \begin{pmatrix}  v_1\\0 \end{pmatrix}\right),
	$$
	$$\left(\begin{pmatrix}   0& \tau_{12}\\\tau_{12}&0 \end{pmatrix}, \begin{pmatrix}  0\\v_2 \end{pmatrix}\right), \quad \left( \begin{pmatrix}   0& 0\\0&\tau_{22}\end{pmatrix}, \begin{pmatrix}  v_1\\0 \end{pmatrix}\right), \quad \left( \begin{pmatrix}   0& 0\\0&\tau_{22} \end{pmatrix}, \begin{pmatrix}  0\\v_2 \end{pmatrix}\right).$$
	
	Finally, we need to solve a least-squares problem with the linear system generated by \eqref{prb11}.
	
	\section{Numerical examples}
	
	In this section, we show the results of RNN-PG and M-RNN-PG methods for solving the linear elasticity problem in two or three dimensions. We compare these methods with FEM, Mixed DG method, and PINNs.

	\begin{example}\label{exam1}In this example, we solve a  2-dimensional linear elasticity problem with a solution 
		$\bu = \begin{pmatrix}
			e^{x-y}xy(1-x)(1-y) \\
			sin(\pi x)sin(\pi y)
		\end{pmatrix}$ 
		and homogeneous Dirichlet boundary condition,
		\begin{numcases}{}
			\bsigma =2\mu \bvarepsilon(\bu) + \lambda \,\rm{tr}\left( \bvarepsilon(\bu)\right) I_2\;\;\;{\rm in} \; \Omega\;  , \notag \\
			- {\rm div}\, \bsigma = \fb  \;\;\;  \quad\;\;\;\; \qquad\quad\;\,{\rm in}  \; \Omega\;  ,\notag \\
			\bu= \bzero\qquad\;\;\;\,\; \qquad \quad\;\;\;\;\;\;\;\;\;{\rm on} \; \Gamma \;,\notag
		\end{numcases}
		on $\Omega = (0,1)^2$ with $\Gamma = \partial \Omega $. Here, $I_2$ is the $2 \times 2$ identity matrix, and the Lamé constants are set to be $\mu = 1/2$ and $\lambda = 1$. 
	\end{example}


	We approximate $\bu$ by a two-layer randomized neural network with $\mathcal{U}(-1, 1)$ as the initial uniform distribution in the RNN-PG method. The network has $2$ input neurons ($n_0=2$) and $2$ output neurons ($n_2=2$), and uses $\rho=\tanh$ as the activation function. In M-RNN-PG methods, we introduce another two-layer randomized neural network to approximate $\bsigma$. This network has 2 input neurons ($n_0=2$) and $3$ output neurons ($n_2=3$). To evaluate the impact of the degrees of freedom (DoF) of our method, we set the number of neurons in the hidden layer to be $n_{1}=100,200,400$, respectively. The DoF for the RNN-PG method is $2n_1$, and the DoF for the M-RNN-PG method is $5n_1$.

	
	\begin{table}[!htbp]	
		\centering  
		\renewcommand{\arraystretch}{1.4}
		\setlength\tabcolsep{1.2mm}
		\begin{tabular}{|c|c|c|c|c|c|c|c|c|c|}  
			\hline  
			\diagbox [width=12em] {Scheme}{$n_1$}&  
			\multicolumn{3}{c|}{100}&\multicolumn{3}{c|}{200}&
			\multicolumn{3}{c|}{400}\cr\cline{1-10}  
			&DoF&$\frac{\Vert \bu-\bu_{\rho} \Vert_0}{\Vert \bu \Vert_0}$ & $\frac{\Vert \bsigma-\bsigma_{\rho}  \Vert_0}{\Vert \bsigma \Vert_0}$ &DoF&$\frac{\Vert \bu-\bu_{\rho} \Vert_0}{\Vert \bu \Vert_0}$& $\frac{\Vert \bsigma-\bsigma_{\rho}  \Vert_0}{\Vert \bsigma \Vert_0}$&DoF&$\frac{\Vert \bu-\bu_{\rho} \Vert_0}{\Vert \bu \Vert_0}$& $\frac{\Vert \bsigma-\bsigma_{\rho}  \Vert_0}{\Vert \bsigma \Vert_0}$\cr\hline  
			RNN-PG method&200&6.367e-8  &3.112e-7&400&7.087e-9 &2.883e-8&800&4.269e-9  &1.489e-8 \cr\hline  
			M-RNN-PG method-1&500&8.448e-8 &2.161e-7&1000&1.024e-8 &2.547e-8&2000&3.372e-9  &1.405e-8\cr\hline  
			M-RNN-PG method-2&500&5.509e-8 &2.300e-7&1000&8.530e-9 &1.413e-8&2000&2.655e-9 &2.400e-9\cr\hline  
			M-RNN-PG method-3&500&5.346e-8&3.817e-7&1000&1.132e-8 &3.913e-8&2000&5.366e-9 	&1.260e-8\cr\hline  
			M-RNN-PG method-4&500&5.552e-8&1.058e-7&1000&1.703e-9&5.661e-9&2000&4.133e-10 	&8.070e-10\cr\hline  
		\end{tabular}  
		\caption{Relative $L^2(\Omega)$ errors of different RNN-PG methods and various $n_1$ in Example \ref{exam1}.}
		\label{table1a}
	\end{table}

	\begin{table}[!htbp] 	
		\centering  
		\renewcommand{\arraystretch}{1.4}
		\setlength\tabcolsep{2.3mm}
		\begin{tabular}{|c|c|c|c|c|c|c|c|c|c|}  
			\hline  
			\diagbox [width=6em] {$h$}{Scheme}&  
			\multicolumn{3}{c|}{$P_1$ linear FEM}&\multicolumn{3}{c|}{$P_2$ quadratic FEM}&
			\multicolumn{3}{c|}{$P_3$ cubic FEM}\cr\cline{1-10}  
			&DoF&$\frac{\Vert \bu-\bu_{h} \Vert_0}{\Vert \bu \Vert_0}$ & $\frac{\Vert \bsigma-\bsigma_{h}  \Vert_0}{\Vert \bsigma \Vert_0}$ &DoF&$\frac{\Vert \bu-\bu_{h} \Vert_0}{\Vert \bu \Vert_0}$& $\frac{\Vert \bsigma-\bsigma_{h}  \Vert_0}{\Vert \bsigma \Vert_0}$&DoF&$\frac{\Vert \bu-\bu_{h} \Vert_0}{\Vert \bu \Vert_0}$& $\frac{\Vert \bsigma-\bsigma_{h}  \Vert_0}{\Vert \bsigma \Vert_0}$ \cr\hline  
			$2^{-2}$&18& 1.815e-1 &3.645e-1&98 &9.035e-3 &5.605e-2&242&7.092e-4 &5.825e-3\cr\hline  
			$2^{-3}$&98&4.643e-2 &1.895e-1&450&1.121e-3 & 1.470e-2&1058&4.095e-5 & 7.367e-4 \cr\hline  
			$2^{-4}$&450&1.184e-2&9.588e-2 &1922&1.386e-4&3.731e-3&4418&2.463e-6 & 9.213e-5 \cr\hline  
			$2^{-5}$&1922&2.980e-3 &4.810e-2&7938& 1.726e-5& 9.366e-4 & 18050&1.514e-7&1.151e-5\cr\hline  
		\end{tabular}  
		\caption{Relative $L^2(\Omega)$ errors of the $P_k$ FEM with different $h$ in Example \ref{exam1}.}
		\label{table1b}
	\end{table}

	The test functions are selected from the bilinear basis functions of the finite element method, which partitions the domain $\Omega$ into square elements with a mesh size of $h=2^{-4}$. The numerical integration is done by the Gauss-Legendre quadrature with $25$ points in each square, and $100$ points are randomly drawn on each edge of $\Gamma$. We calculate $\nabla \bu_\rho$ by the central difference method with a spacing of $10^{-6}$, which achieves a good balance between computation time and accuracy. We use the least-square solver $scipy.linalg.lstsq$ in Python to solve the linear system that arises.


	We show the relative $L^2$ errors for $\bu$ and $\bsigma$ with different values of $n_1$ and different versions of RNN-PG methods in Table \ref{table1a}. We find that the accuracy improves as the DoF increases. In addition, M-RNN-PG method-4 attains a very low relative $L^2$ error of around $10^{-10}$, due to its benefit of not using numerical derivatives and imposing boundary conditions. For comparison, we also calculate the numerical solutions $\bu_h$ and  $\bsigma_h$ by FEM with Fenics (\cite{Langtangena2017fenics}) on a sequence of uniform triangulations with mesh size $h=2^{-n}\; (n=2,3,4,5)$, and we employ the standard triangle Lagrange elements $P_k$ ($k=1,2,3$) to solve this problem. Here, $k$ is the polynomial degree, and the DoF for $P_k$ FEM is $2(k/h-1)^2$. The numerical errors are given in Table \ref{table1b}. By comparing Table \ref{table1a} and Table \ref{table1b}, we notice that both the RNN-PG method and M-RNN-PG methods surpass FEM in terms of DoF and accuracy, which implies that our method can achieve a more accurate numerical representation with much fewer degrees of freedom.

\begin{figure}[!htbp] 		
		\centering
		\includegraphics[scale=0.6]{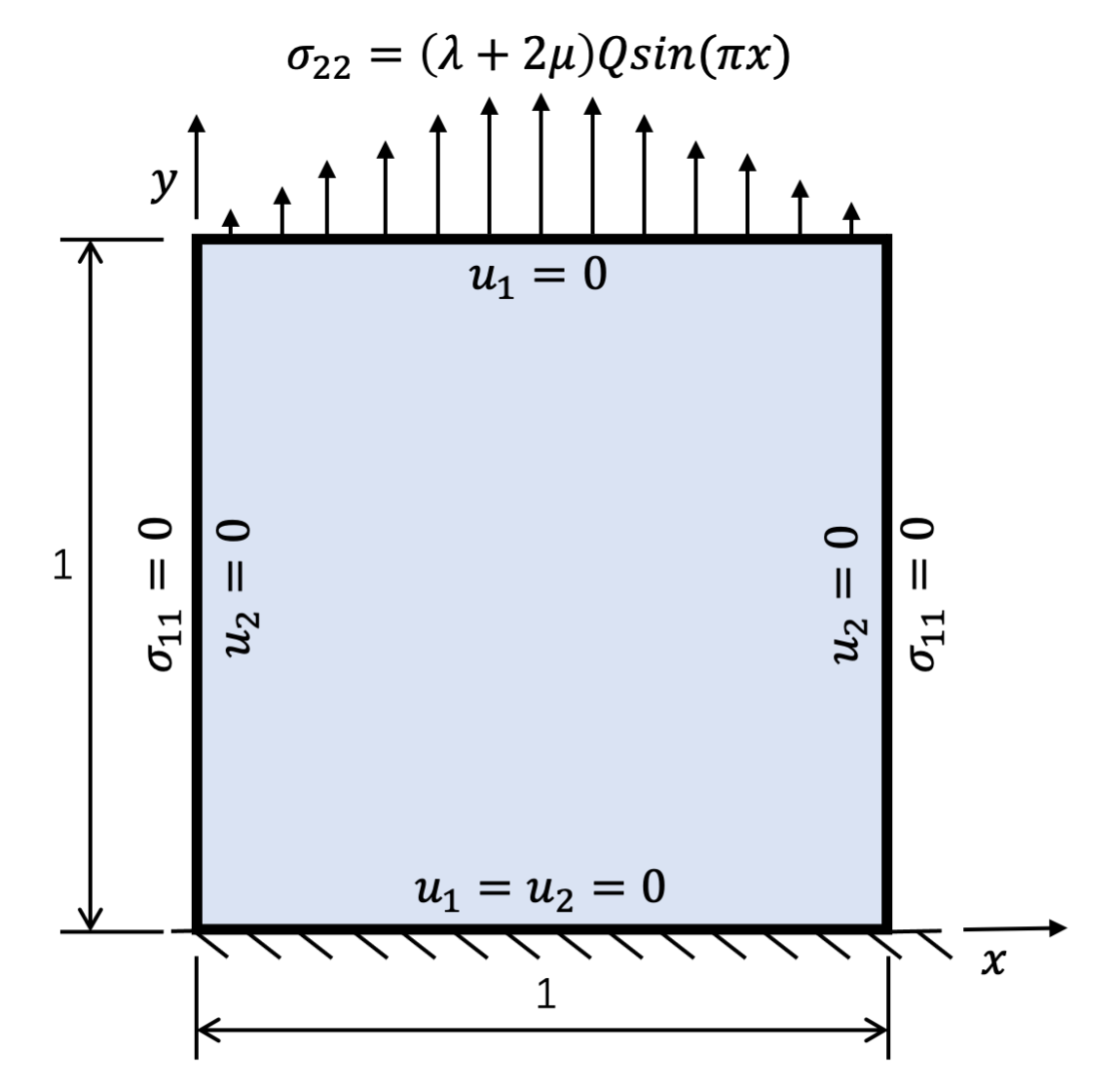}
		\caption{Problem setup and boundary conditions in Example \ref{exam2}. }
		\label{setup}
	\end{figure}

	\begin{example}\label{exam2}
		We consider a  2-dimensional elasticity plane-strain problem with an exact solution	$\bu = \begin{pmatrix}
			cos(2\pi x)sin(\pi y) \\
			sin(\pi x)Qy^4/4
		\end{pmatrix}$  on the unit square, and the boundary conditions are depicted in the {\rm Figure \ref{setup}}.
		\begin{numcases}{}
			\bsigma =2\mu \bvarepsilon(\bu) + \lambda \,\rm{tr} \left( \bvarepsilon(\bu)\right) I_2\;\;\;\; {\rm in} \; \Omega\;  , \notag \\
			- {\rm div}\, \bsigma = \fb  \;\;\;  \quad\;\;\;\; \qquad\quad\;\;\,{\rm in}  \; \Omega\;  ,\notag 
		\end{numcases}
		where $I_2$ is the $2 \times 2$ identity matrix, and parameters are set to be $\mu = 1/2$ , $\lambda = 1$ and $Q = 4$. 
	\end{example}
	
\begin{figure}[!ht]
		\begin{center}
			\subfigure[Numerical solution $u_{\rho1}$]{
				\centering
				\includegraphics[width=2.0in]{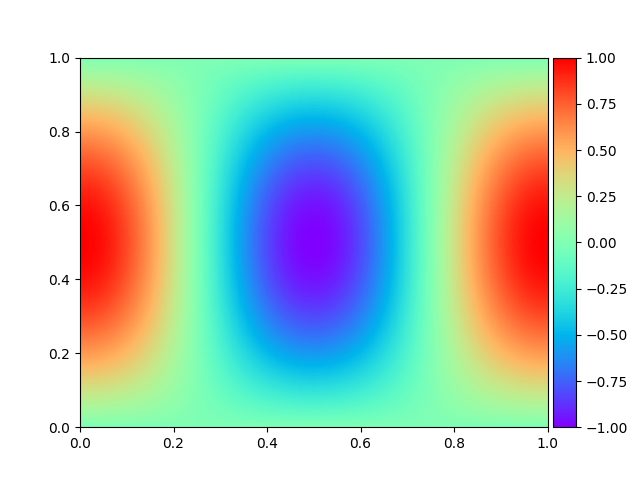}
			}
			\subfigure[Exact solution $u_1$]{
				\centering
				\includegraphics[width=2.0in]{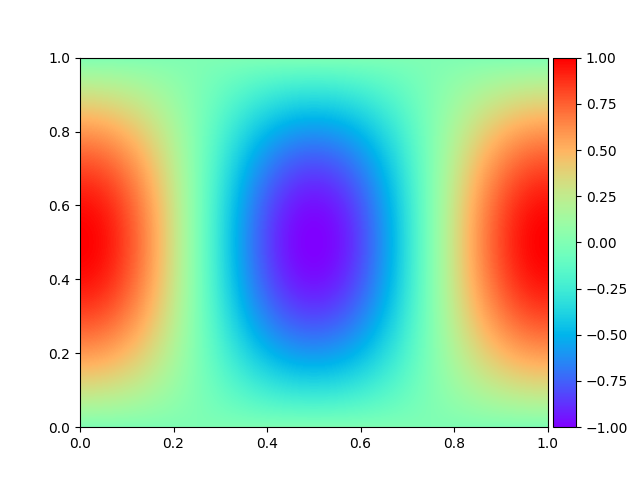}
			}
			\subfigure[$|u_1-u_{\rho1}|$]{
				\centering
				\includegraphics[width=2.0in]{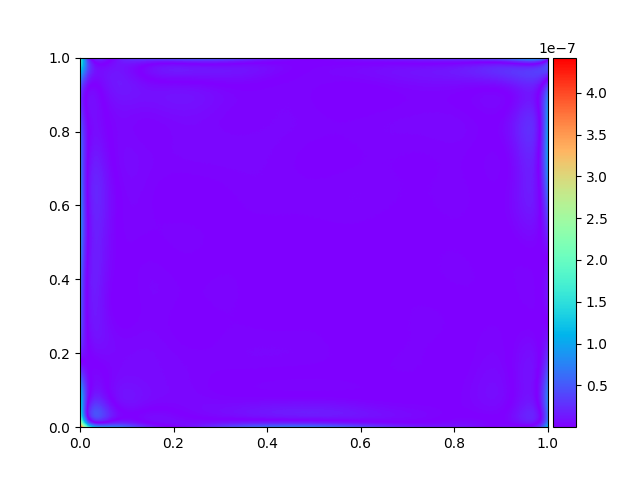}
			}
			
			\subfigure[Numerical solution $u_{\rho2}$]{
				\centering
				\includegraphics[width=2.0in]{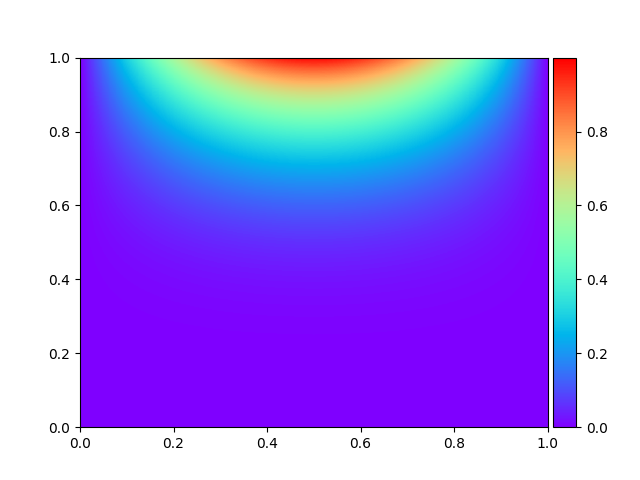}
			}
			\subfigure[Exact solution $u_2$]{
				\centering
				\includegraphics[width=2.0in]{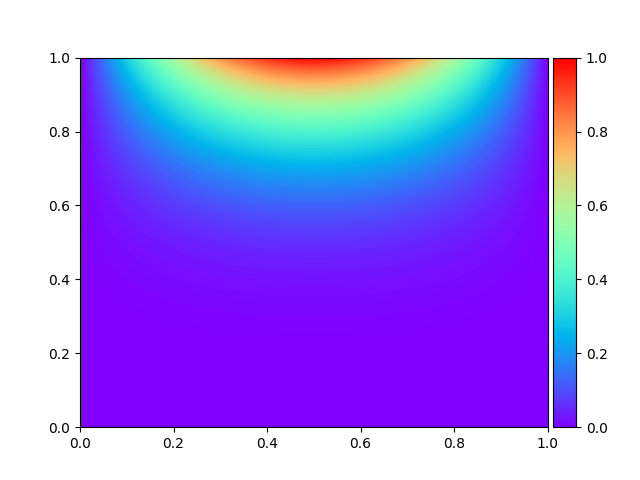}
			}
			\subfigure[$|u_2-u_{\rho2}|$]{
				\centering
				\includegraphics[width=2.0in]{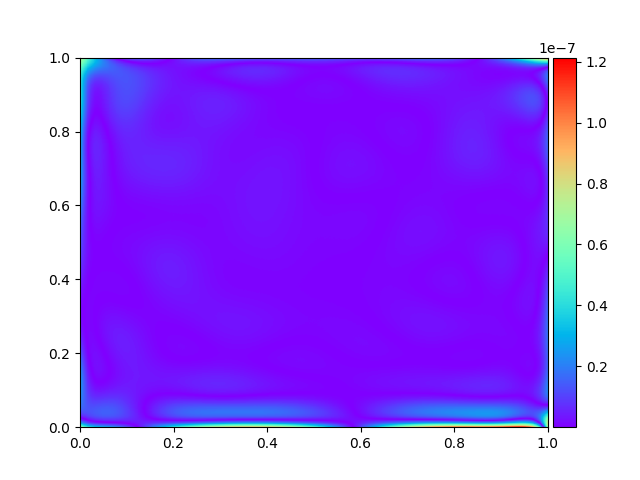}
			}
			
		\end{center}
		\vspace*{-15pt}
		\caption{Numerical solution $\bu_{\rho}$ by M-RNN-PG method-4 with $h=2^{-4}$ and  DoF $=2000$ in Example \ref{exam2}.}
		\label{figure2}
	\end{figure}

\begin{figure}[!ht]
		\begin{center}
			\subfigure[Numerical solution $\sigma_{\rho11}$]{
				\centering
				\includegraphics[width=2.0in]{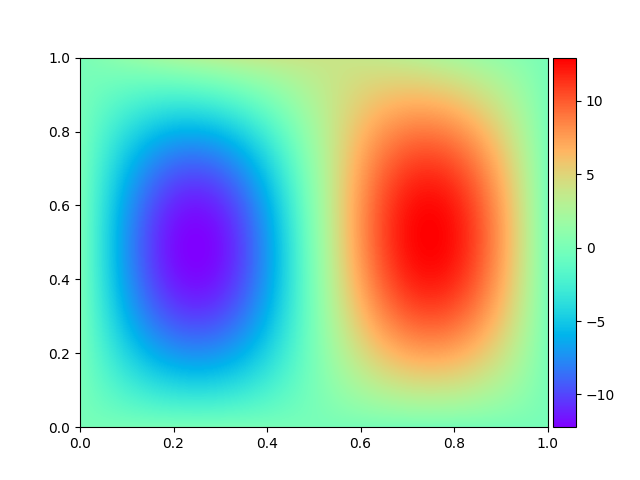}
			}
			\subfigure[Exact solution $\sigma_{11}$]{
				\centering
				\includegraphics[width=2.0in]{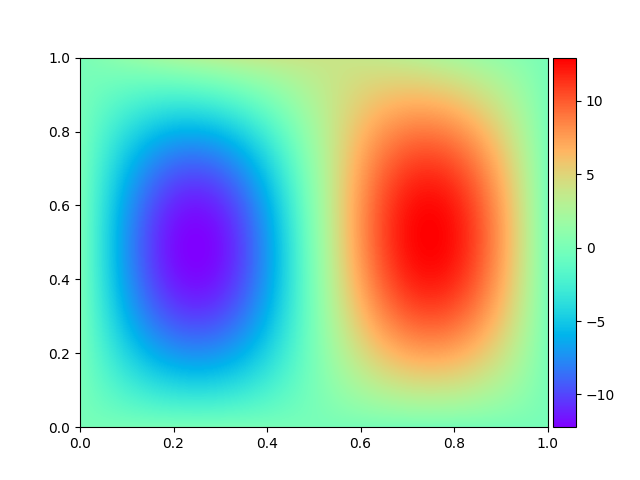}
			}
			\subfigure[$|\sigma_{11}-\sigma_{\rho11}|$]{
				\centering
				\includegraphics[width=2.0in]{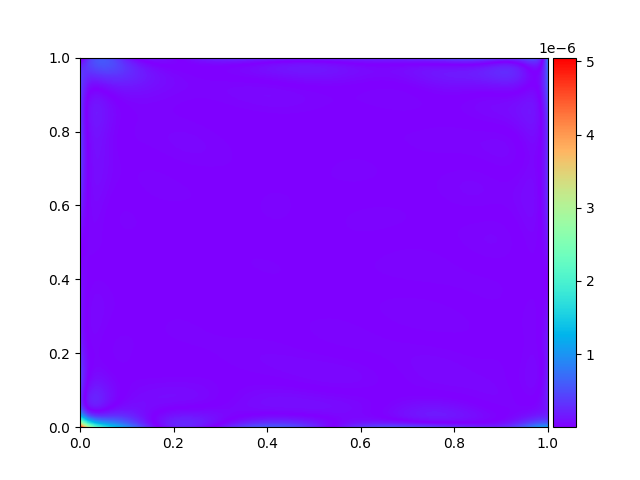}
			}
			
			\subfigure[Numerical solution $\sigma_{\rho12}$]{
				\centering
				\includegraphics[width=2.0in]{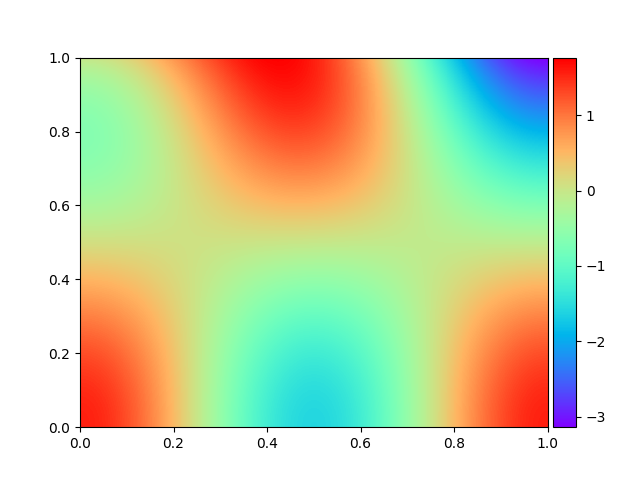}
			}
			\subfigure[Exact solution $\sigma_{12}$]{
				\centering
				\includegraphics[width=2.0in]{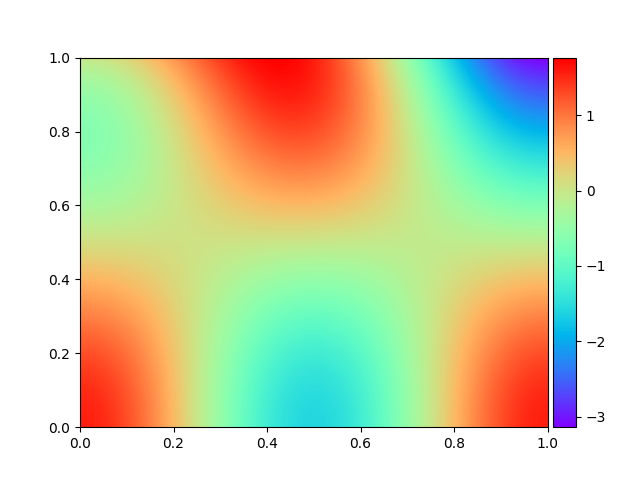}
			}
			\subfigure[$|\sigma_{12}-\sigma_{\rho12}|$]{
				\centering
				\includegraphics[width=2.0in]{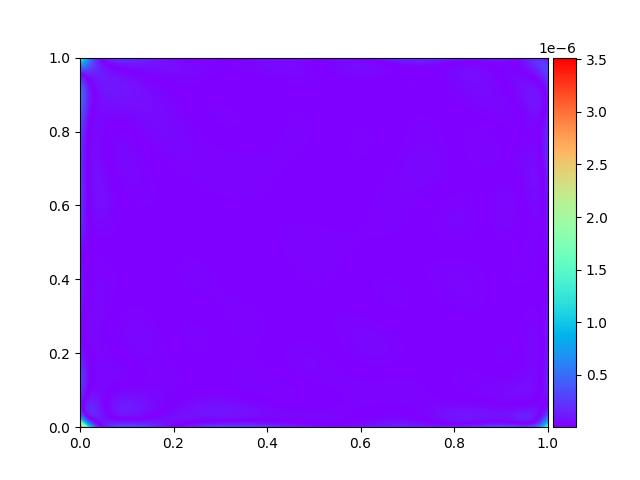}
			}
			
			\subfigure[Numerical solution $\sigma_{\rho22}$]{
				\centering
				\includegraphics[width=2.0in]{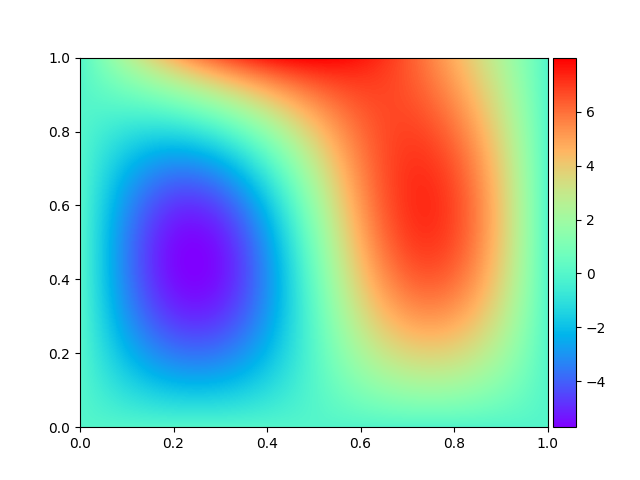}
			}
			\subfigure[Exact solution $\sigma_{22}$]{
				\centering
				\includegraphics[width=2.0in]{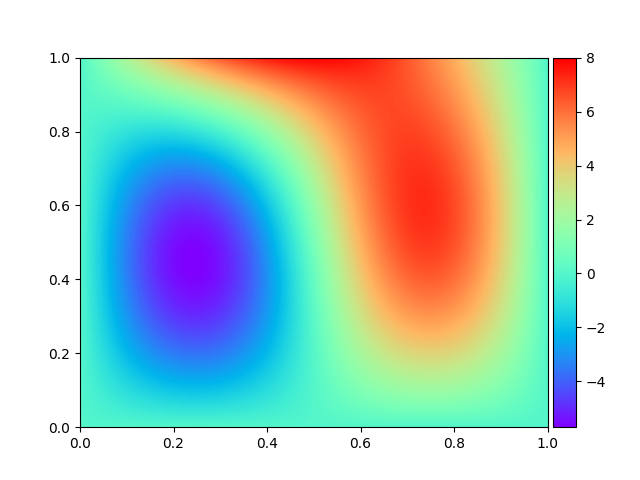}
			}
			\subfigure[$|\sigma_{22}-\sigma_{\rho22}|$]{
				\centering
				\includegraphics[width=2.0in]{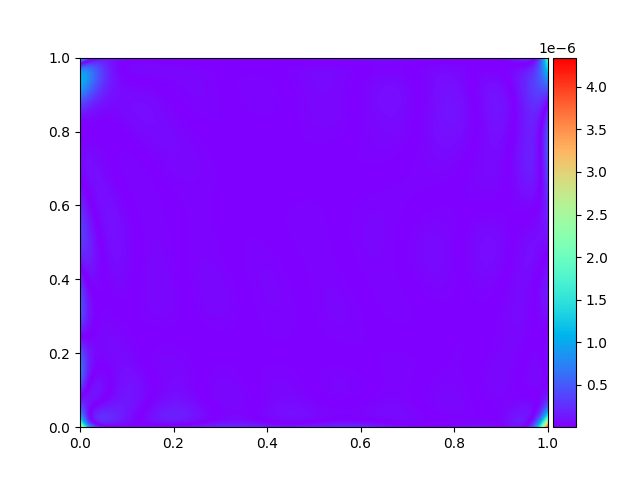}
			}
			
		\end{center}
		\vspace*{-15pt}
		\caption{Numerical solution $\bsigma_{\rho}$ by M-RNN-PG method-4 with $h=2^{-4}$ and  DoF $=2000$  in Example \ref{exam2}.}
		\label{figure3}
	\end{figure}

	We use the same network structure and test functions as in Example \ref{exam1}, and we report the relative errors in Table \ref{table2}. Compared to the results in \cite{Guo2022ex2elasticity,Haghighat2020ex2solid} obtained by PINNs, which require a lot of time to train deep neural networks to achieve the relative $L^2$ errors around $10^{-3}$, both RNN-PG method and M-RNN-PG methods can reach even $10^{-8}$ with much less time by using randomized neural networks. Meanwhile, we notice that all the variants of RNN-PG methods can handle this complex boundary condition very well, and we display the numerical solution of $\bu_\rho  = (  u_{\rho1},u_{\rho2})^t$ and $\bsigma_\rho = ( \sigma_{\rho11},\sigma_{\rho12},\sigma_{\rho22})^t$ obtained by M-RNN-PG methods-4 in Figure \ref{figure2} and Figure \ref{figure3}, respectively.
	
	
		\begin{table}[!htbp]	
		\centering  
		\renewcommand{\arraystretch}{1.4}
		\setlength\tabcolsep{1.2mm}
		\begin{tabular}{|c|c|c|c|c|c|c|c|c|c|}  
			\hline  
			\diagbox [width=12em] {Scheme}{$n_1$}&  
			\multicolumn{3}{c|}{100}&\multicolumn{3}{c|}{200}&
			\multicolumn{3}{c|}{400}\cr\cline{1-10}  
			&DoF&$\frac{\Vert \bu-\bu_{\rho} \Vert_0}{\Vert \bu \Vert_0}$ & $\frac{\Vert \bsigma-\bsigma_{\rho}  \Vert_0}{\Vert \bsigma \Vert_0}$ &DoF&$\frac{\Vert \bu-\bu_{\rho} \Vert_0}{\Vert \bu \Vert_0}$& $\frac{\Vert \bsigma-\bsigma_{\rho}  \Vert_0}{\Vert \bsigma \Vert_0}$&DoF&$\frac{\Vert \bu-\bu_{\rho} \Vert_0}{\Vert \bu \Vert_0}$& $\frac{\Vert \bsigma-\bsigma_{\rho}  \Vert_0}{\Vert \bsigma \Vert_0}$\cr\hline  
			RNN-PG method&200&6.623e-6  &1.167e-5&400&3.818e-7  &7.979e-7&800&9.422e-8  &3.065e-7  \cr\hline  
			M-RNN-PG method-1&500&5.912e-6  &6.202e-6&1000&8.130e-7 &1.287e-6&2000&2.998e-7  &3.155e-7\cr\hline  
			M-RNN-PG method-2&500&6.486e-6 &1.418e-5&1000&3.153e-6  &1.020e-6&2000&7.949e-7  &7.746e-7\cr\hline  
			M-RNN-PG method-3&500&6.550e-6&6.433e-6&1000&3.504e-6 &3.678e-6&2000&7.161e-7 	&4.833e-7\cr\hline  
			M-RNN-PG method-4&500&8.795e-7&1.570e-6&1000&7.559e-8 &1.300e-7&2000&2.318e-8&2.145e-8\cr\hline  
		\end{tabular}  
		\caption{Relative $L^2(\Omega)$ errors of different RNN-PG methods and various $n_1$ in Example \ref{exam2}.}
		\label{table2}
	\end{table}

	
	
	\begin{example}\label{exam3}
		
		In this example, we consider a two-dimensional problem of linear elasticity
		\begin{numcases}{}
			\bsigma =2\mu \bvarepsilon(\bu) + \lambda \,\rm{tr} \left( \bvarepsilon(\bu)\right) I_2\;\;\; {\rm in} \; \Omega\;  , \notag \\
			- {\rm div}\, \bsigma = \fb  \;\;\;  \quad\;\;\;\; \qquad\quad\;\,{\rm in}  \; \Omega\;  ,\notag \\
			\bu= \bzero\qquad\;\;\;\,\; \qquad \quad\;\;\;\;\;\;\;\;\;{\rm on} \; \Gamma \;,\notag
		\end{numcases}
		on domain $\Omega = (0,1)^2$ with a solution 
		$\bu = \begin{pmatrix}
			-x^2(x-1)^2y(y-1)(2y-1) \\
			x(x-1)(2x-1)y^2(y-1)^2.
		\end{pmatrix}$ 
		Here, $I_2$ is the identity matrix of size $2 \times 2$, and $\mu$ and $\lambda$ are the Lam\'e parameters. 
		These parameters depend on Young’s modulus $E$ and the Poisson’s ratio $\nu$ of the material, as follows:
		$$\mu = \frac{E}{2(1+\nu)},\qquad
		\lambda = \frac{E\nu}{(1+\nu)(1-2\nu)}.$$
		The Poisson’s ratio $\nu$ measures the compressibility of the material, and it approaches 0.5 as the material becomes more incompressible. 
	\end{example}
	
	In this example, we use the same network structure as in Example \ref{exam1}. We choose test functions on a square mesh with mesh size $h=2^{-5}$ to increase the amount of data, and we set the spacing parameter of the central difference method to $10^{-8}$. 
	We show the $L^2$ errors for $\bu$ and $\bsigma$ with various values of $\nu$, DoF by mixed RNN-PG methods in Table \ref{table3a}. We also compare our results with those of the mixed DG method (\cite{Wang2020mixdg}) in Table \ref{table3b}, where we change Poisson’s ratio $\nu$, $h$, and the elements of the method. We observe that our method is locking-free in the incompressible limit case, similar to the mixed DG method, highly accurate but with fewer DoF than Mixed DG method. In addition, we see that M-RNN-PG method-4 has better accuracy for $\bsigma$ than the other methods when larger $\nu$ is considered.

	\begin{table}[!htbp] 	
		\centering  
		\renewcommand{\arraystretch}{1.4}
		\setlength\tabcolsep{0.5mm}
		\begin{tabular}{|c|c|c|c|c|c|c|c|c|c|}  
			\hline  
			\multicolumn{2}{|c|}{\diagbox [width=9em] {$\nu$\;,\; DoF}{Scheme}}
			&  
			\multicolumn{2}{m{3cm}<{\centering}|}{{M-RNN-PG method 1}} 
			&	\multicolumn{2}{m{3cm}<{\centering}|}{{M-RNN-PG method 2}} &
			\multicolumn{2}{m{3cm}<{\centering}|}{{M-RNN-PG method 3}} & 	\multicolumn{2}{m{3cm}<{\centering}|}{{M-RNN-PG method 4}} \cr\cline{1-10}  
			$\nu$&DoF&$\Vert \bu-\bu_{\rho} \Vert_0$ & $\Vert \bsigma-\bsigma_{\rho} \Vert_0$ &$\Vert \bu-\bu_{\rho} \Vert_0$ & $\Vert \bsigma-\bsigma_{\rho} \Vert_0$  &
			$\Vert \bu-\bu_{\rho} \Vert_0$ & $\Vert \bsigma-\bsigma_{\rho} \Vert_0$  &
			$\Vert \bu-\bu_{\rho} \Vert_0$ & $\Vert \bsigma-\bsigma_{\rho} \Vert_0$\cr\hline  
			\multirow{4}{*}{$ \nu = 0.49$}&$500$&2.611e-7 &6.904e-5 &7.840e-8&6.370e-6 &1.655e-7&	1.142e-5&5.881e-9	&4.153e-7 	\cr\cline{2-10}
			&$1000$& 4.288e-8&	1.421e-6&1.363e-8&	1.448e-6&4.216e-8 &	2.454e-6&2.872e-10	&3.076e-8    
			\cr\cline{2-10}
			&$2000$& 5.494e-8	&1.015e-5&5.164e-9 	&4.966e-7&4.599e-8 	&7.026e-6 &1.478e-10	&1.415e-8    
			\cr\cline{2-10}
			&$4000$& 3.109e-8	&4.485e-6&9.060e-9	&8.345e-7&2.253e-8	&2.181e-6&5.942e-11 	&5.791e-9     
			\cr\hline  
			\multirow{4}{*}{$ \nu = 0.4999$}&$500$&1.504e-7	&4.202e-3&1.156e-7	&1.193e-4&3.238e-7	&5.701e-3&4.796e-9 	&2.712e-6     
			\cr\cline{2-10} 
			&$1000$& 4.726e-8	&6.491e-4&1.133e-8 	&8.830e-6&3.522e-8	&6.223e-4&2.307e-10	&2.637e-7 
			\cr\cline{2-10}
			&$2000$& 3.763e-8 &	6.772e-4&4.483e-9&	4.755e-6&2.615e-8	&7.912e-4&1.464e-10 &	4.444e-7
			\cr\cline{2-10}
			&$4000$& 4.594e-8	&9.960e-4 &8.364e-9	&1.764e-6 &6.228e-8	&1.515e-3&1.067e-10	&4.216e-8
			\cr\hline  
			\multirow{4}{*}{$ \nu = 0.499999$}&$500$& 2.524e-7	&2.644e-1&6.464e-8	&4.242e-3&2.158e-7	&2.776e-2 &3.854e-9	&4.248e-4    
			\cr\cline{2-10}
			&$1000$&9.581e-8 &	1.149e-1 &2.318e-8	&3.092e-3&3.582e-8&	2.092e-2&2.484e-10	&7.245e-6     
			\cr\cline{2-10}
			&$2000$&3.946e-8	&8.892e-2&7.034e-9 	&3.337e-4&4.623e-8	&3.840e-3 &1.086e-10	&8.319e-6     
			\cr\cline{2-10}
			&$4000$& 2.811e-8 	&3.430e-2&2.060e-8	&1.144e-5&4.940e-8	&7.410e-2&9.871e-11	&3.203e-6    
			\cr\hline  
			
		\end{tabular}  
		\caption{$L^2(\Omega)$ error of different formulations of RNN-PG methods with different $\nu$ and DoF in Example \ref{exam3}.}
		\label{table3a}
	\end{table}
	
	\begin{table}[!htbp] 	
		\centering  
		\renewcommand{\arraystretch}{1.4}
		\setlength\tabcolsep{0.5mm}
		\begin{tabular}{|c|c|c|c|c|c|c|c|c|c|c|}  
			\hline  
			\multicolumn{2}{|c|}{\diagbox [width=10em] {$\nu$\;, \;$h$}{Scheme}}&  
			\multicolumn{3}{c|}{$P_1^{-1}-P_0^{-1}$ }&\multicolumn{3}{c|}{$P_2^{-1}-P_1^{-1}$} & \multicolumn{3}{c|}{$P_3^{-1}-P_2^{-1}$} \cr\cline{1-11}  
			$\nu$&$h$ &DoF	&$\Vert \bu-\bu_{h} \Vert_0$ & $\Vert \bsigma-\bsigma_{h} \Vert_0$ &DoF &$\Vert \bu-\bu_{h} \Vert_0$ & $\Vert \bsigma-\bsigma_{h} \Vert_0$  & DoF &$\Vert \bu-\bu_{h} \Vert_0$ & $\Vert \bsigma-\bsigma_{h} \Vert_0$  \cr\hline  		
			\multirow{2}{*}{$\nu = 0.49$}& $2^{-4}$ & 5632& 4.196e-4 &1.649e-3 &12288 &2.940e-5&6.685e-5&21504&1.282e-6&1.674e-6 \cr\cline{2-11}		
			&$2^{-5}$&22528&2.106e-4&7.906e-4&49152&7.380e-6&1.434e-5&86016&1.610e-7&1.110e-7 \cr\hline  	
			\multirow{2}{*}{$\nu = 0.4999$}& $2^{-4}$ & 5632 & 4.196e-4&1.640e-3&12288 &2.940e-5&6.652e-5&21504&1.282e-6&1.665e-6 \cr\cline{2-11}		
			& $2^{-5}$&22528 & 2.106e-4&7.864e-4&49152 &7.380e-6&1.427e-5&86016&1.610e-7&1.110e-7\cr\hline  
			\multirow{2}{*}{$\nu = 0.499999$}&  $2^{-4}$ & 5632 &4.196e-4&1.640e-3&12288 &2.940e-5&6.652e-5&21504&1.282e-6&1.665e-6 \cr\cline{2-11}
			& $2^{-5}$& 22528 &2.106e-4&7.863e-4&49152 &7.380e-6&1.427e-5&86016&1.610e-7&1.110e-7\cr\hline  
			
		\end{tabular}  
		\caption{$L^2(\Omega)$ error of the mixed DG method (taken from \cite{Wang2020mixdg}) with different $\nu$ and $h$ in Example \ref{exam3}.}
		\label{table3b}
	\end{table}

	\begin{example}\label{exam4}we consider a three-dimensional linear elasticity problem
		\begin{numcases}{}
			\bsigma =2\mu \bvarepsilon(\bu) + \lambda \,\rm{tr} \left( \bvarepsilon(\bu)\right) I_3\;\;\; \,{\rm in} \; \Omega\;  , \notag \\
			- {\rm div}\, \bsigma = \fb  \;\;\;  \quad\;\;\;\; \qquad\quad\;\,{\rm in}  \; \Omega\;  ,\notag \\
			\bu= \bg\qquad\;\;\;\,\; \qquad \quad\;\;\;\;\;\;\;\;\;{\rm on} \; \Gamma_D \;,\notag
		\end{numcases}
		on the unit cube with an exact solution	
		$$\bu = \begin{pmatrix}
			2^4 \\
			2^5\\
			2^6
		\end{pmatrix}x(1-x)y(1-y)z(1-z).
		$$ 
		Here, $I_3$ is the $3 \times 3$ identity matrix, and the Lam\'e  parameters are set to be $\mu = 1/2$ ,$\lambda = 1$. 
	\end{example}

	In this example, we employ a two-layer randomized neural network with 3 input neurons ($n_0 = 3$) and 3 output neurons ($n_2 = 3$) to approximate $\bu$, which is initialized with a uniform distribution $\mathcal{U}(-1, 1)$. Our test functions are trilinear functions on cubic meshes with mesh size $h=2^{-n}\; (n=2,3,4)$. We perform the numerical integration by using Gauss-Legendre quadrature with 125 points inside each cube. We randomly sample 100 points on each face of $\Gamma_D$ to impose the Dirichlet boundary condition. Table \ref{table4a} shows the $L^2$ errors for $\bu$ and $\bsigma$ with different mesh sizes and numbers of DoF. We calculate $\bsigma_{\rho}$ by the relation $\bsigma_{\rho} =2\mu \bvarepsilon(\bu_{\rho}) + \lambda \,\rm{tr} \left( \bvarepsilon(\bu_{\rho})\right) I_3$. The RNN-PG method solves this three-dimensional problem very well and achieves high accuracy with fewer degrees of freedom compared to the results obtained by FEM and mixed DG method (\cite{Wang2020mixdg}) in Table \ref{table4b}.

	
	\begin{table}[!htbp]	
		\centering  
		\renewcommand{\arraystretch}{1.4}
		\setlength\tabcolsep{0.4mm}
		\begin{tabular}{|c|c|c|c|c|c|c|c|c|}  
			\hline  
			\diagbox [width=6em] {$h$}{DoF}&  \multicolumn{2}{c|}{300}&
			\multicolumn{2}{c|}{600}&\multicolumn{2}{c|}{1200}&
			\multicolumn{2}{c|}{2400}\cr\cline{1-9}  
			&$\Vert \bu-\bu_{\rho} \Vert_0$ & $\Vert \bsigma-\bsigma_{\rho} \Vert_0$  
			&$\Vert \bu-\bu_{\rho} \Vert_0$ & $\Vert \bsigma-\bsigma_{\rho} \Vert_0$ 
			&$\Vert \bu-\bu_{\rho} \Vert_0$ & $\Vert \bsigma-\bsigma_{\rho} \Vert_0$ 
			&$\Vert \bu-\bu_{\rho} \Vert_0$ & $\Vert \bsigma-\bsigma_{\rho} \Vert_0$ \cr\hline  
			$2^{-2}$&1.125e-2&1.444e-1&1.266e-3&4.172e-2&1.878e-3&5.460e-2&1.538e-3&3.547e-2 \cr\hline  
			$2^{-3}$&1.732e-2&2.123e-1&5.036e-4&6.410e-3&5.086e-6&8.350e-5&1.210e-6 &1.740e-5\cr\hline  
			$2^{-4}$&7.251e-2&7.193e-1&7.902e-4&8.954e-3&5.909e-6&1.019e-4&5.519e-8&1.207e-6\cr\hline  
		\end{tabular}  
		\caption{$L^2(\Omega)$  error of the RNN-PG method with different $h$ and DoF in Example \ref{exam4}.}
		\label{table4a}
	\end{table}
	
			%
	
	\begin{table}[!htbp] 	
		\centering  
		\renewcommand{\arraystretch}{1.4}
		\setlength\tabcolsep{3mm}
		\begin{tabular}{|c|c|c|c|c|c|c|c|c|}  
			\hline  
			\diagbox [width=6em] {$h$}{Scheme}&  
			\multicolumn{3}{c|}{$P_2$  FEM}
			&\multicolumn{3}{c|}{$P_2^{-1}-P_1^{-1}$}  \cr\cline{1-7}
			&DoF&$\Vert \bu-\bu_{h} \Vert_0$ & $\Vert \bsigma-\bsigma_{h} \Vert_0$ &DoF& $\Vert \bu-\bu_{h} \Vert_0$ & $\Vert \bsigma-\bsigma_{h} \Vert_0$  \cr\hline  
			$2^{-1}$&81&4.858e-2&1.085&3456&8.310e-2&3.642e-1\cr\hline  
			$2^{-2}$&1029&6.511e-3&3.401e-1&27648 &2.274e-2&6.646e-2\cr\hline  
			$2^{-3}$&10125&8.232e-4&8.359e-2& 221184 &5.820e-3&1.238e-2\cr\hline  
			
		\end{tabular}  
		\caption{$L^2(\Omega)$ errors of the FEM and mixed DG method  (taken from \cite{Wang2020mixdg}) with different $h$ in Example \ref{exam4}.}
		\label{table4b}
	\end{table}

	\section{Summary}

	In this paper, we explore the use of RNN-PG methods for solving linear elasticity problems in solid mechanics. Based on the weak formulation, we can naturally incorporate the PDE information and the Neumann boundary condition, and we can enforce the Dirichlet boundary condition by using random samples. We apply the least-squares method to solve the resulting linear system, and we obtain a more accurate solution of the displacement and the stress tensor. Moreover, we use M-RNN-PG methods that employ separate neural networks to approximate displacement and stress variables, and we can easily achieve the symmetric property of the stress tensor in the network structure. Numerical results show that traditional methods such as FEM and DG method require a larger degree of freedom to obtain an accurate numerical solution, while our method provides a more accurate solution with less DoF. Compared to the common PINNs, both methods can handle problems with complex boundary conditions, complex domains, or higher dimensions, but our method has the advantage of being more accurate and less time-consuming, as solving a linear system is much faster than training a deep neural network.
	
	Despite the success of RNN-PG methods and M-RNN-PG methods, several questions remain, such as how to generate a proper initialization with a suitable activation function, how to compute the resulting linear system efficiently, especially in higher dimensions with a large condition number, and how to extend the proposed approach to other mechanical problems such as nonlinear elasticity, viscoplasticity, and elastoplasticity.

	\bibliographystyle{siam}
	
\end{document}